\documentclass[review]{elsarticle}

\journal{a journal}

\usepackage{lineno,hyperref}
\modulolinenumbers[5]
\usepackage{graphicx}
\usepackage{epstopdf}
\usepackage[utf8]{inputenc}
\usepackage{amsmath}
\usepackage{amstext}
\usepackage{amsfonts}
\usepackage{amsthm}
\usepackage{natbib,amssymb}
\usepackage[bf,SL,BF]{subfigure}
\usepackage{caption}
\usepackage{color}
\usepackage{listings}

\hypersetup{colorlinks}

\newcommand{\RR}{{\mathbb{R}}}
\newcommand{\EE}{{\mathbb{E}}}
\newcommand{\PP}{{\mathbb{P}}} 
\newcommand{\F}{{\mathcal{F}}} 

\newtheorem{definition}{Definition}[section]
\newtheorem{theorem}[definition]{Theorem}
\newtheorem{lemma}[definition]{Lemma}

\newtheorem{assumption}[definition]{Assumption}

\newcommand{\eproof}{\hfill $\Box$}  

\begin{document}

\begin{frontmatter}

\title{Strong convergence rate of Euler-Maruyama method for stochastic differential equations with H\"older continuous drift coefficient driven by symmetric $\alpha$-stable process}

\author{Wei Liu}
\ead{weiliu@shnu.edu.cn; lwbvb@hotmail.com}
\address{Department of Mathematics, Shanghai Normal University, Shanghai, 200234, China}

\begin{abstract}
Euler-Maruyama method is studied to approximate stochastic differential equations driven by the symmetric $\alpha$-stable additive noise with the $\beta$ H\"older continuous drift coefficient. When $\alpha \in (1,2)$ and $\beta \in (0,\alpha/2)$, for $p \in (0,2]$ the $L^p$ strong convergence rate is proved to be $\beta/\alpha$. The proofs in this paper are extensively based on H\"older's and Bihari's inequalities, which is significantly different from those in Huang and Liao (2018).
\end{abstract}


\end{frontmatter}


\section{Introduction} \label{sec:intro}
In this paper, we consider the Euler-Maruyama (EM) method for $d$-dimensional stochastic differential equations (SDEs) driven by the symmetric $\alpha$-stable process
\begin{equation}\label{eq:SDE}
dx(t) = f(x(t))dt + d L(t),~~~x(0) = x_0 \in \RR^d,
\end{equation}
where the drift coefficient $f: \RR^d \rightarrow \RR^d$ is $\beta$ H\"older continuous and $L(t)$ is a scalar symmetric $\alpha$-stable process. Throughout this paper, we assume that $\alpha \in (1,2)$. When $\alpha = 2$, $L(t)$ is the standard Brownian motion.
\par
In Chapter 1 of \cite{ST1994book}, the authors present four equivalent ways to describe the $\alpha$-stable process. In this paper, we adopt the following description that
\begin{itemize}
\item $L(0) = 0$ a.s.;
\item $L(t)$ has independent increments;
\item $L(t) - L(s)$ follows $S_\alpha((t-s)^{1/\alpha},0,0)$ for any $0 \leq s < t < \infty$, where $S_\alpha(\sigma, \beta, \mu)$ is a four-parameter stable distribution.
\end{itemize}
It should be mentioned that such a description makes numerical simulations quite straightforward (see Section \ref{sec:NumExp} for more details). The symmetric $\alpha$-stable process belongs to the family of L\'evy processes. We refer the readers to \cite{Applebaum2009book} for the detailed introduction to L\'evy processes driven SDEs.
\par
Since the explicit expressions of the true solutions are hardly obtained, numerical methods become extremely important. When the driven noise is the standard Brownian motion, numerical methods for SDEs under the standard assumptions on coefficients are well studied \cite{KP1992a}. When the driven noise is $\alpha$-stable process, authors in \cite{JMW1996} investigated the EM method under the standard assumptions on coefficients. More recently, $\alpha$-stable process driven SDEs with the H\"older's continuous drift have been attracting a lot of attention. The existence and uniqueness of \eqref{eq:SDE} was studied in \cite{Priola2012}. When the driven noise $L(t)$ is the truncated symmetric $\alpha$-stable process, the strong convergence rate of the EM method was given in \cite{PT2017}. When the driven noise $L(t)$ is the symmetric $\alpha$-stable process the strong convergence rate of the EM method was proved in \cite{HL2018}. The proofs in both \cite{PT2017} and \cite{HL2018} are based on the associated Kolmogorov equation.
\par
In this paper, we present an alternative proof of the strong convergence for \eqref{eq:SDE}, which extensively uses inequalities, such as H\"older's inequality and Bihari's inequality.
\par
The main differences between our paper and \cite{HL2018} are as follows.
\begin{itemize}
\item Our approach works on the difference between the true and numerical solutions directly without the knowledge of the associated Kolmogorov equation.
\item We do not need the drift coefficient to be bounded.
\end{itemize}
To keep the notation simple and to present our ideas clearly, we only investigate the case of the additive scalar noise in this paper. However, the techniques used in this paper can be extended to the case of the multiplicative multi-dimensional noise as well as the case of asymmetric $\alpha$ stable processes. Due to the limit of pages, we will report relevant results in further works.
\par
This paper is constructed as follows. The assumptions and the main result are presented in Section \ref{sec:main}. The proofs are given in Section \ref{sec:proofs}. A numerical example is displayed in Section \ref{sec:NumExp}.

\section{Assumptions and the main result} \label{sec:main}
Throughout this paper, unless otherwise specified, we let $(\Omega,\F,\PP)$ be a complete probability space with a filtration $\{\F_t\}_{t\in [0,T]}$ satisfying the usual conditions (that is, it is right continuous and increasing while $\F_0$ contains all $\PP$-null sets), and let $\EE$ denote the probability expectation with respect to $\PP$. If $x\in \RR^d$, then $|x|$ is the Euclidean norm. Moreover, for two real numbers $a$ and $b$, we use $a\vee b=max(a,b)$.

\begin{assumption}\label{ass:HolderDrift}
Assume that there exist constants $K > 0$ and $\beta \in (0,1)$ such that
\begin{equation*}
|f(x) - f(y)| \leq K |x - y|^\beta,
\end{equation*}
for any $x,y\in \RR^d$.
\end{assumption}
From Assumption \ref{ass:HolderDrift}, it is easy to see that
\begin{equation} \label{cond:lineargrowth}
|f(x)| \leq K_2 (1 + |x|^\beta),
\end{equation}
where $K_2 = K \vee |f(0)|$.

For some given time step $\Delta \in (0,1)$ and the terminal time $T$, define $N = T/\Delta$. The EM method to the SDE (\ref{eq:SDE}) is defined by
\begin{equation} \label{eq:disEMmethod}
Y_{i+1} = Y_i + f(Y_i) \Delta + \Delta L_i,
\end{equation}
where $Y_0 = x_0$ and $\Delta L_i$ follows the stable distribution $S_\alpha(\Delta^{1/\alpha},0,0)$ for $i=1,2,...,N$ \cite{JMW1996}. Here, $Y_i$ is the approximation to $x(i \Delta)$, for $i=1,2,...,N$. We also define the continuous version of \eqref{eq:disEMmethod} by
\begin{equation}\label{eq:continuousEMmethod}
Y(t) = Y(0) + \int_0^t f(\bar{Y}(s))ds + L(t), ~~~Y(0) = Y_0 = x_0,
\end{equation}
where $\bar{Y}(t) = Y_i$, when $t \in [i\Delta, (i+1)\Delta)$.
\par
Our main result is as follows.
\begin{theorem}\label{theorem:strongrate}
Suppose that Assuption \ref{ass:HolderDrift} holds. If $2 \beta < \alpha$, for any $p \in (0,2]$ the strong error of the EM method \eqref{eq:continuousEMmethod} is
\begin{equation*}
\sup_{0 \leq t \leq T} \EE |x(t) - Y(t)|^p \leq \Delta^{p\beta/\alpha} C_5^{p/2} \exp(C_6Tp/2),
\end{equation*}
where $C_5 = 2 K^2 T^{3/2}C_4^{2\beta/q}$ and $C_6 = 2 K^2 T^{1/2}\left(\left( C_5 \Delta^{2\beta/\alpha} \right)^{1-\beta} + (1 - \beta) 2 K^2 T^{3/2} \right)^{-1}$.
\end{theorem}

\section{Proofs} \label{sec:proofs}
To prove Theorem \ref{theorem:strongrate}, we first present three lemmas.
\begin{lemma}\label{lemma:SDEqthmomentbd}
Suppose that Assumption \ref{ass:HolderDrift} holds. For any $q \in [1, \alpha)$, the $q$th moment of the solution to \eqref{eq:SDE} is bounded
\begin{equation*}
\sup_{0 \leq t \leq T} \EE |x(t)|^q \leq C_3,
\end{equation*}
where
\begin{equation*}
C_3 = \left[C_2^{1-\beta} + (1 - \beta) \left(6K_2 \right)^q T\right]^{1/(1-\beta)},
\end{equation*}
\begin{equation*}
C_2 = 3^q ( \EE |x(0)|^q + (2 K_2)^q T^{(2q - 1)/q} + C_1^q T^{q/\alpha} ),
\end{equation*}
and $C_1$ is a constant dependent on $\alpha$ and $q$ (see Property 1.2.17 at Page 18 of \cite{ST1994book} for the exact expression of $C_1$).
\end{lemma}
\begin{proof}
For any $q \in [1,\alpha)$ and any $t \in [0,T]$, by the elementary inequality we derive from \eqref{eq:SDE} that
\begin{equation*}
|x(t)|^q \leq 3^q \left( |x(0)|^q + \left| \int_0^t f(x(s)) ds \right|^q  + |L(t)|^q \right).
\end{equation*}
By H\"older's inequality (see for example page 5 in \cite{M2008a}), we have
\begin{equation*}
|x(t)|^q \leq 3^q \left( |x(0)|^q +  T^{(q-1)/q}\int_0^t \left|f(x(s))\right|^q ds   + |L(t)|^q \right).
\end{equation*}
Taking expectations both sides, by \eqref{cond:lineargrowth} and the elementary inequality we obtain
\begin{equation*}
\EE |x(t)|^q \leq 3^q \left( \EE |x(0)|^q + (2K_2)^{q} T^{(q-1)/q} \int_0^t \left( 1 + \EE |x(s)|^{q \beta} \right) ds + C_1 T^{q/\alpha}  \right),
\end{equation*}
where the fact $\EE |L(t)| ^q \leq C_1 t^{q/\alpha}$ (see Property 1.2.17 at Page 18 of \cite{ST1994book}) is used.
\par
Applying H\"older's inequality for $\EE |x(s)|^{q \beta}$, we can see
\begin{equation*}
\EE |x(t)|^q \leq \left( 6 K_2 \right)^q \int_0^t \left(  \EE |x(s)|^{q} \right)^{\beta} ds + C_2,
\end{equation*}
where $C_2 = 3^q ( \EE |x(0)|^q + (2 K_2)^q T^{(2q - 1)/q} + C_1^q T^{q/\alpha} )$.
\par
By Bihari's inequality (see for example Page 45 in \cite{M2008a}), the assertion holds.
\end{proof}

\begin{lemma}\label{lemma:EMsolqthmomentbd}
Suppose that Assumption \ref{ass:HolderDrift} holds. For any $q \in [1, \alpha)$, the $q$th moment of the solution to the EM method \eqref{eq:continuousEMmethod} is bounded
\begin{equation*}
\sup_{0 \leq t \leq T} \EE |Y(t)|^q \leq C_3,
\end{equation*}
where $C_3$ is the same as that in Lemma \ref{lemma:SDEqthmomentbd}.
\end{lemma}
\begin{proof}
For any $q \in [1, \alpha)$ and any $t \in [0,T]$, following the similar steps as those in the proof of Lemma \ref{lemma:SDEqthmomentbd}, we can get
\begin{equation*}
\EE |Y(t)|^q \leq \left( 6 K_2 \right)^q \int_0^t \left(  \EE |\bar{Y}(s)|^{q} \right)^{\beta} ds + C_2.
\end{equation*}
Then taking  the supremum inside the integral on the right hand side, we can see
\begin{equation*}
\EE |Y(t)|^q \leq \left( 6 K_2 \right)^q \int_0^t \left( \sup_{0 \leq u \leq s} \EE |Y(u)|^{q} \right)^{\beta} ds + C_2.
\end{equation*}
Since the inequality above holds for any $t \in [0,T]$, we have
\begin{equation*}
\sup_{0\leq v \leq t}\EE |Y(v)|^q \leq \left( 6 K_2 \right)^q \int_0^t \left( \sup_{0 \leq u \leq s} \EE |Y(u)|^{q} \right)^{\beta} ds + C_2.
\end{equation*}
Applying Bihari's inequality, the proof is completed.
\end{proof}

\begin{lemma}\label{lemma:Y-Ybar}
Suppose that Assumption \ref{ass:HolderDrift} holds. For any $q \in [1, \alpha)$ and any $t \in [0,T]$, the difference between the continuous and discrete versions of the EM method is
\begin{equation*}
\EE |Y(t) - \bar{Y} (t)|^q  \leq C_4 \Delta^{q/\alpha},
\end{equation*}
where $C_4 =  2^{2q}  K_2^q (1 + C_3^\beta) + 2^q C_1^q$.
\end{lemma}

\begin{proof}
From \eqref{eq:continuousEMmethod}, for any $t \in [i\Delta, (i+1)\Delta)$ we have
\begin{equation*}
Y(t) - \bar{Y} (t) = \int_{i\Delta}^t f(\bar{Y}(s))ds + (L(t)-L(i\Delta)).
\end{equation*}
For any $q \in [1, \alpha)$, taking the power of $q$ on both sides we can get, in the similar manner as Lemma \ref{lemma:SDEqthmomentbd}, that
\begin{equation}\label{pf:Y-Yb}
|Y(t) - \bar{Y} (t)|^q \leq 2^{2q} \Delta^{(2q-1)/q} K_2^q (1 + |\bar{Y}(t)|^{q\beta}) + 2^q  |L(t)-L(i\Delta)|^q.
\end{equation}
Due to the selfsimilarity of the symmetric $\alpha$-stable process, we have
\begin{equation*}
\EE |L(t)-L(i\Delta)|^q = \EE |L(t- i\Delta)|^q \leq C_1^q \Delta^{q/\alpha}.
\end{equation*}
Taking expectations on both sides of \eqref{pf:Y-Yb} yields
\begin{equation*}
\EE |Y(t) - \bar{Y} (t)|^q \leq \Delta^{q/\alpha} \left( 2^{2q}  K_2^q (1 + \EE |\bar{Y}(t)|^{q\beta}) + 2^q C_1^q  \right),
\end{equation*}
where $(2q-1)/q > q/\alpha$ is used. Since $q\beta \in (0,q)$, by H\"older's inequality and Lemma \ref{lemma:EMsolqthmomentbd} we have
\begin{equation*}
\sup_{0 \leq t \leq T}\EE |\bar{Y}(t)|^{q\beta} \leq \sup_{0 \leq t \leq T} \left( \EE |\bar{Y}(t)|^q \right)^{\beta} \leq \sup_{0 \leq t \leq T} \left( \EE |Y(t)|^q \right)^{\beta} \leq C_3^\beta.
\end{equation*}
Therefore, the assertion holds.
\end{proof}

Now we are ready to prove the main result.

{\it \textbf{Proof of Theorem \ref{theorem:strongrate}}}
\par \noindent
From \eqref{eq:SDE} and \eqref{eq:continuousEMmethod}, we have
\begin{equation*}
x(t) - Y(t) = \int_0^t \left( f(x(s)) - f(Y(s)) \right)ds + \int_0^t \left( f(Y(s)) - f(\bar{Y}(s)) \right)ds.
\end{equation*}
Taking squares on both sides, by Assumption \ref{ass:HolderDrift} and H\"older's inequality we have
\begin{equation}\label{pf:x-Y}
|x(t) - Y(t)|^2 \leq 2 K^2 T^{1/2} \left( \int_0^t |x(s) - Y(s)|^{2\beta} ds +  \int_0^t |Y(s) - \bar{Y}(s)|^{2\beta} ds \right),
\end{equation}
for any $t \in [0,T]$. Since $2 \beta < \alpha$, we can choose a $q$ such that $q \in ( 2\beta, \alpha)$. Then by H\"older's inequality and Lemma \ref{lemma:Y-Ybar}, we can see
\begin{equation*}
\EE |Y(t) - \bar{Y}(t)|^{2 \beta} \leq \left( \EE |Y(t) - \bar{Y}(t)|^q \right)^{2\beta/q} \leq C_4^{2\beta/q} \Delta^{2\beta/\alpha},
\end{equation*}
and
\begin{equation*}
\EE|x(t) - Y(t)|^{2\beta} \leq \left( \EE |x(t) - Y(t)|^2 \right)^{\beta}.
\end{equation*}
Thus, taking expectations on both sides of \eqref{pf:x-Y} we have
\begin{equation}\label{pf:Ex-y}
\EE |x(t) - Y(t)|^2 \leq 2 K^2 T^{1/2} \int_0^t \left( \EE |x(s) - Y(s)|^2 \right)^{\beta} ds + C_5 \Delta^{2\beta/\alpha},
\end{equation}
where $C_5 = 2 K^2 T^{3/2}C_4^{2\beta/q}$. Using Bihari's inequality, we obtain
\begin{equation} \label{pf:bdonEx-Y}
\EE |x(t) - Y(t)|^2 \leq  \left(\left( C_5 \Delta^{2\beta/\alpha} \right)^{1-\beta} + (1 - \beta) 2 K^2 T^{3/2} \right)^{1/(1-\beta)}.
\end{equation}
Now we rewrite the right hand side of \eqref{pf:Ex-y} into
\begin{eqnarray*}
&&\EE |x(t) - Y(t)|^2 \nonumber \\
&\leq& 2 K^2 T^{1/2} \int_0^t \EE |x(s) - Y(s)|^2 \times \left( \EE |x(s) - Y(s)|^2 \right)^{\beta-1} ds + C_5 \Delta^{2\beta/\alpha}.
\end{eqnarray*}
By \eqref{pf:bdonEx-Y}, we can see
\begin{equation*}
\EE |x(t) - Y(t)|^2 \leq C_6 \int_0^t \EE |x(s) - Y(s)|^2 ds + C_5 \Delta^{2\beta/\alpha},
\end{equation*}
where $C_6 = 2 K^2 T^{1/2}\left(\left( C_5 \Delta^{2\beta/\alpha} \right)^{1-\beta} + (1 - \beta) 2 K^2 T^{3/2} \right)^{-1}$. By Gronwall's inequality (see for example Page 45 in \cite{M2008a}), we have
\begin{equation*}
\EE |x(t) - Y(t)|^2 \leq \Delta^{2\beta/\alpha} C_5 \exp(C_6 T).
\end{equation*}
For any $p \in (0,2]$, applying H\"older's inequality results in the assertion.  \eproof

\section{Numerical example} \label{sec:NumExp}
To make the EM method for \eqref{eq:SDE} implementable to those readers who are interested in computer simulations, we use the scalar SDE
\begin{equation}\label{eq:example}
dx(t) = x^{4/9}(t)dt + L(t),~~~x(0) = 1,
\end{equation}
as an example.
\par
\begin{figure}[htbp]
\centering
\subfigure[Probability density functions]
{
  \begin{minipage}{5cm}
  \label{fig:left}
  \centering
  \includegraphics[scale=0.25]{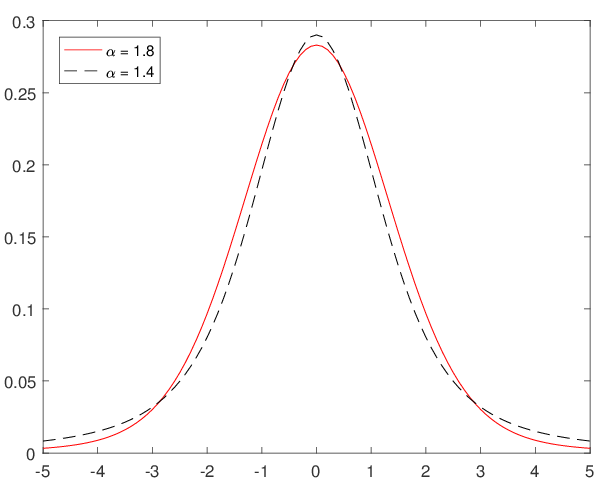}
  \end{minipage}
}
\subfigure[Two paths of the EM solutions]
{
  \begin{minipage}{5cm}
  \label{fig:right}
  \centering
  \includegraphics[scale=0.25]{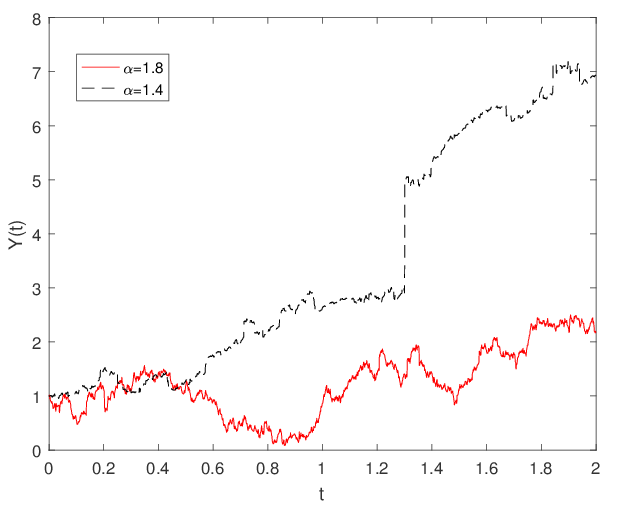}
  \end{minipage}
}
\caption{Probability density functions of the driven noise and paths of EM solutions}
\end{figure}

The matlab codes to generate one path of the EM solution to \eqref{eq:example} with the time step $\Delta = 0.001$, $\alpha = 1.8$ and $T=2$ are as follows.
\begin{lstlisting}[frame=single]
T=2;   h=0.001;  % terminal time T and step size h
N=T/h;           % number of iterations
a=1.8;           % value of alpha
sp=h^(1/a);  % value of the scale parameter
dL=random('stable',a,0,sp,0,1,N); % generate noise
X=zeros(1,N+1); % vector to contain the solution
X(1)=1;   % initial value
for i = 1:N
  X(i+1) = X(i) + h*X(i)^(4/9) + dL(i); % EM method
end
\end{lstlisting}
Figure \ref{fig:left} shows the probability density functions of the symmetric stable distribution $S_\alpha(1,0,0)$ with $\alpha = 1.8$ and $\alpha = 1.4$, respectively. It can be seen that when $\alpha$ gets smaller the tails become heavier, which means higher probability is allocated to values far away from the centre. Two sample paths of \eqref{eq:example} with $\alpha = 1.8$ and $\alpha = 1.4$ are plotted in Figure \ref{fig:right}. It can be observed that the path with $\alpha = 1.4$ has larger jumps than the path with $\alpha = 1.8$, which is due to the heavier tails of the distribution.

\section*{Acknowledgement}
Wei Liu is financially supported by the National Natural Science Foundation of China (11701378, 11871343) and “Chenguang Program” supported by both Shanghai Education Development Foundation and Shanghai Municipal Education Commission (16CG50).

\end{document}